\documentclass[journal,twoside,web]{ieeecolor}

\usepackage{generic}
\usepackage{cite}
\usepackage{amsmath,amssymb,amsfonts}
\usepackage{textcomp}

\usepackage{amsmath,amssymb,amsfonts}
\usepackage{algcompatible}
\usepackage{algorithm}
\usepackage{cite}

\usepackage{bm}
\usepackage{multirow}

\usepackage{graphicx}
\usepackage{subfigure}
\usepackage{epsfig} 
\usepackage{cite}
\usepackage{lcsys}

\newtheorem{theorem}{Theorem}
\newtheorem{lemma}{Lemma}

\newtheorem{definition}[theorem]{Definition}
\newtheorem{remark}{Remark}
\newtheorem{assumption}{Assumption}

\renewcommand{\QED}{$\blacksquare$}
\newcommand{\norm}[1]{\left\lVert#1\right\rVert}

\newcommand{\xx}{\mathbf{x}}
\newcommand{\bb}{\mathbf{b}}

\newcommand{\zz}{\mathbf{z}}
\renewcommand{\AA}{\mathbf{A}}
\newcommand{\FF}{\mathbf{F}}
\newcommand{\FFt}{\mathbf{F}_t}
\newcommand{\FFo}{\overline{\FF}_t}
\newcommand{\vv}{\bm{\nu}}
\newcommand{\vt}{\mathbf{v}_t}

\newcommand{\DD}{\mathbf{D}}

\newcommand{\xh}{\mathbf{\hat{x}}}

\begin{document}

\def\BibTeX{{\rm B\kern-.05em{\sc i\kern-.025em b}\kern-.08em
    T\kern-.1667em\lower.7ex\hbox{E}\kern-.125emX}}
\markboth{\journalname, VOL. XX, NO. XX, XXXX 2017}
{Author \MakeLowercase{\textit{et al.}}: Preparation of Papers for IEEE Control Systems Letters (August 2022)}

\title{An Online Newton's Method for Time-varying Linear Equality Constraints}

\author{Jean-Luc Lupien and
Antoine Lesage-Landry, IEEE Member
\thanks{J.-L. Lupien and A. Lesage-Landry are with the Department of Electrical Engineering, Polytechnique Montréal, GERAD \& Mila, Montréal, QC, Canada, H3T 1J4. e-mail: \{\texttt{jean-luc.lupien, antoine.lesage-landry}\}\texttt{@polymtl.ca}.}
\thanks{This work was funded by the Institute for Data Valorization (IVADO) and the Natural Sciences and Engineering Research Council (NSERC).}%
}

\maketitle
\thispagestyle{empty}

\begin{abstract}
We consider online optimization problems with time-varying linear equality constraints. In this framework, an agent makes sequential decisions using only prior information. At every round, the agent suffers an environment-determined loss and must satisfy time-varying constraints. Both the loss functions and the constraints can be chosen adversarially. We propose the Online Projected Equality-constrained Newton Method (OPEN-M) to tackle this family of problems. We obtain sublinear dynamic regret and constraint violation bounds for OPEN-M under mild conditions. Namely, smoothness of the loss function and boundedness of the inverse Hessian at the optimum are required, but not convexity. Finally, we show OPEN-M outperforms state-of-the-art online constrained optimization algorithms in a numerical network flow application.\end{abstract}


\section{Introduction}
\IEEEPARstart{I}{n} online convex optimization (OCO), an agent aims to sequentially play the best decision with respect to a potentially adversarial loss function using only prior information \cite{Shalev,zinkevich2003online}. {In other words, decisions must be made before observing the loss function and constraints.} OCO algorithms have many applications including portfolio selection, artificial intelligence, and real-time control of power systems \cite{taylor2016power,hazanOCO,applications}. 

The preponderant performance metric for OCO algorithms is regret \cite{Shalev}, the cumulative difference between the loss incurred by the agent and that of a comparator sequence. Two main types of regret exist: static and dynamic. For static regret, the comparator sequence is defined as the best fixed decision in hindsight \cite{zinkevich2003online}. For dynamic regret, the comparator sequence is the round-optimal decision in hindsight \cite{zinkevich2003online, hazanOCO}. In OCO, one intends to design an algorithm that possesses a sublinear regret bound. Sublinear regret implies that the time-averaged regret goes to zero as time increases. The OCO algorithm therefore plays the best decision with respect to the comparator sequence over sufficiently long time horizons \cite{Shalev, hazanOCO}.

A persistent obstacle in OCO algorithmic design has been the integration of time-varying constraints. In this context, the decision sequence must also satisfy environment-determined constraints \cite{MOSP}. The performance of such algorithms is measured both in terms of its regret but also its constraint violation, the cumulative distance from feasibility of the agent's decisions. Similarly to the regret, a sublinear constraint violation is desired and implies that, on average, decisions will be feasible over a large time horizon \cite{JMLR}. When considering time-varying constraints, an analysis using the dynamic regret is preferable over one based on the static regret because the best fixed feasible decision can have an arbitrarily large loss or might not even exist \cite{FTRLTime-Varying}. 

Most OCO algorithms tackling time-varying constraints are analysed based only on static regret. In \cite{JMLR,FTRLTime-Varying} and \cite{HazanInterior}, sublinear regret and violation bounds are achieved for long-term constraints. Sublinear static regret and constraint violation bounds are also achieved in \cite{OnlineEpsilonProof} using virtual queues and in \cite{OnlineConvexBandit} using an online saddle-point algorithm. More recently, an augmented Lagrangian method \cite{LagrangianTime-Varying} has been shown to outperform previous Lagragian-based methods \cite{MOSP, OnlineConvexBandit, OnlineEpsilonProof} in numerical experiments.

Algorithms with dynamic regret bounds have also been developed such as the modified online saddle-point method (\texttt{MOSP}) \cite{MOSP}. This algorithm has simultaneous sublinear dynamic regret and constraint violation bounds. However, \texttt{MOSP}'s bounds are dependent on strict conditions on optimal primal and dual variable variations in addition to time-sensitive step sizes. An exact-penalty method for dealing with time-varying constraints possessing sublinear regret and constraint violation bounds was presented in \cite{L1penalty} but shares \texttt{MOSP}'s step size limitation. Virtual queues are used in \cite{QueueConstrained} with time-varying constraints achieving simultaneous sublinear regret and constraint violation bounds without requiring Slater's condition to hold. Sublinear regret and constraint violation are achieved in \cite{DynamicPolyak} including for some non-convex functions but with considerably looser bounds.

The application of an interior-point method to time-varying convex optimization is presented in \cite{InteriorOnline}. However, this context differs from OCO because current-round information is available to the decision-maker.

{In this work, our main contribution is the design of a novel online optimization algorithm that can efficiently deal with time-varying linear equality constraints in the OCO setting. This method simultaneously possesses the tightest dynamic regret and constraint violation bounds presented thus far for constrained OCO problems. Additionally, the method does not require hyperparameters, time-dependent step sizes, or a predefined time horizon making it easily implementable.}

\section{Background}

In recent work, a second-order method for online optimization yielded tighter dynamic regret bounds compared to first-order approaches \cite{OnlineNewton}. Specifically, \cite{OnlineNewton} proposes an online extension of Newtons' method, \texttt{ONM}, applicable to non-convex problems that possesses a tight dynamic regret bound. However, this approach is only applicable to unconstrained problems. In an offline setting, the unconstrained and linear equality-constrained Newton's method's performance are the same~\cite{boyd, renegar}. This result motivates the extension of \texttt{ONM} to a setting with time-varying linear equality constraints.

\subsection{Problem definition}

We consider online optimization problems of the following form. Let $\xx_t\in \mathbb{R}^n$, $n\in\mathbb{N}$ be the decision vector at time $t$. Let $f_t$: $\mathbb{R}^n\mapsto \mathbb{R}$ be a twice-differentiable function. Let $\mathbf{A}_t \in \mathbb{R}^{p\times n}$ be a rank $p\in\mathbb{N}$ { full row-rank matrix}, and $\bb_t \in \mathbb{R}^p$. The problem at round $t=1,2,...,T$ can then be written as:
\begin{align}\label{BaseProblem}
\begin{split}
    &\min\limits_{\xx_t} f_t(\xx_t)\\
    \text{s.t.}\quad& \AA_t\xx_t=\bb_t.
\end{split}
\end{align}

In this work, dynamic regret will be used as the performance metric because it is more stringent compared to static regret. Indeed, sublinear dynamic regret implies sublinear static regret \cite{FTRLTime-Varying}. Dynamic regret $R_{\text{d}}(T)$ is defined as:
\begin{align}\label{RegretDef}R_{\text{d}}(T) = \sum\limits_{t=1}^{T} \big[f_t(\xx_t)-f_t(\xx^*_t)\big],\end{align}
where $\xx_t^*$ is the round-optimal solution and $T\in\mathbb{N}$ is the time horizon. For \eqref{BaseProblem}, the round-optimum $\xx_t^*$ is the solution to the following system of equations:
\begin{align*}
    \nabla f_t(\xx^*_t) + \AA_t^\top\vv^*_t&=0\\
    \AA_t\xx^*_t-\bb_t &= 0,
\end{align*}
where $\vv_t^*\in \mathbb{R}^n$ is the dual variable associated with \eqref{BaseProblem}'s equality constraints.

The constraint violation term is defined as:
\begin{align}\label{VioDef}\text{Vio}(T) = \sum\limits_{t=1}^T \norm{\AA_t\xx_t-\bb_t},\end{align}
and quantifies the cumulative distance from feasibility, with respect to the Euclidean norm, of the decision sequence. All norms, $\norm{\cdot}$, refer to the Euclidean norm in the sequel. Constraint violation is zero if decisions are feasible at all rounds. {This definition is similar to that used in \cite{MOSP, DynamicPolyak} and is stricter than in \cite{LagrangianTime-Varying} and \cite{QueueConstrained} because the constraints must be satisfied at every timestep and not on average.}

In the first part of this work, we investigate the case for which the feasible space is the same for all rounds, i.e., $\AA_t=\AA, \bb_t=\bb$ for all $t$. The second part builds on this result and extends it to time-varying equality constraints.

\subsection{Preliminaries}

We define the function $\DD_t(\xx):\mathbb{R}^n\mapsto\mathbb{R}^{(n+p)^2}$ as:
\begin{equation*}\DD_t(\xx) = \begin{bmatrix}\nabla^2f_t(\xx) & \mathbf{A}^\top_t\\\mathbf{A}_t&0\end{bmatrix},\end{equation*}
where $\nabla^2f_t(\xx)$ is the Hessian matrix of $f_t$. We assume that the Hessian is invertible for all $t$ which implies that $\DD_t(\xx)$ is also invertible~\cite[Section 10.1]{boyd}. This guarantees that the Newton update is defined at every round.

Next, we present the online equality-constrained Newton (\texttt{OEN}) update. For any feasible point $\xx_t$, the \texttt{OEN} update minimizes the second-order approximation of $f_t$ around $\xx_t$ subject to the equality constraints. An estimate of the optimal dual variable, $\vv_t$, is also obtained from the update.

\begin{definition}[\texttt{OEN} update]\label{nstepdef}
    The \texttt{OEN} update is:
\begin{align}
\begin{split}
    \begin{bmatrix}\Delta\xx_t\\ \vv_t\end{bmatrix} &= -\DD_t^{-1}(\xx_t)\begin{bmatrix}\nabla f_t(\xx_t)\\0\end{bmatrix}\\
    \xx_{t+1} &= \xx_t+\Delta\xx_t.
    \end{split}
\end{align}
\end{definition}

Let $\vt\in\mathbb{R}^n$ be the difference between subsequent optima:
$\vt =\xx_{t+1}^*-\xx_t^*$.
Throughout this work, we assume:
$0\le\norm{\vt}\le \overline{v}$ for all $t$. {This limits the variation in optima between two subsequent rounds. It is a common assumption in dynamic OCO \cite{MOSP, OnlineNewton, DynamicPolyak}. This assumption could be satisfied in real-world applications such as electric grids where the temporal continuity imposed by the underlying physics limits the variation in optima provided the timestep is sufficiently small.}
The total variation $V_T$ is defined as:
\[V_T = \sum\limits_{t=0}^{T-1} \norm{\xx_{t+1}^*-\xx_t^*}=\sum\limits_{t=0}^{T-1}\norm{\vt},\]
and is bounded above by $V_T\le \overline{v}T$.

An important tool for the analysis of the \texttt{OEN} update is the reduced function $\tilde{f}_t(\zz): \mathbb{R}^{n-p}\mapsto \mathbb{R}$ which is a representation of $f_t(\xx)$ over the feasible set. 

Let $\FF_t \in\mathbb{R}^{n\times(n-p)}$ be such that $\mathcal{R}(\FF_t) = \mathcal{N}(\AA_t)$ where $\mathcal{R}$ is the column space of a matrix and $\mathcal{N}$ is its null space.
Let $\xh\in\mathbb{R}^n$ be such that $\AA_t\xh-\bb_t=\mathbf{0}$. Then, $\tilde{f}_t$ is defined as:
\[\tilde{f_t}(\zz) = f_t(\FF_t \zz + \xh).\]
We remark that the reduced function shares minima with $f_t$, i.e., $\min\limits_{\zz} \tilde{f_t}(\zz) = f_t(\xx_t^*)$ \cite{boyd}. This gives rise to an equivalent unconstrained, reduced problem:
$\min\limits_{\zz} \tilde{f}_t(\zz)$, which can be solved using \texttt{ONM} \cite{OnlineNewton}.

Additionally, there exists a unitary matrix $\overline{\FF}_t\in\mathbb{R}^{n\times(n-p)}$ such that $\mathcal{R}\left(\overline{\FF}_t\right) = \mathcal{N}(\AA_t)$. Without loss of generality, we let $\FF_t=\overline{\FF}_t$ in our analysis. It follows that for any $\zz_t$ satisfying $\xx_t = \overline{\FF}_t\zz_t+\xh$,
we have $\norm{\xx_t-\xx_t^*}=\norm{\overline{\FF}_t(\zz_t-\zz_t^*)}=\norm{\zz_t-\zz_t^*}$.
In other words, the norm of the original problem's \texttt{OEN} update and the reduced problem's \texttt{ONM} update coincide when $\FF_t=\overline{\FF}_t$. 

\subsection{Assumptions}
We now introduce three recurring assumptions that are used throughout this work. These assumptions are mild and notably do not require the objective function to be convex. {These assumptions must hold for all $t=1,2,3,...,T$.} 
\begin{assumption}
\label{ass:inv}
There exists a constant $h>0$ such that: \[\norm{\nabla^2f_t(\xx_t^*)^{-1}}\le\frac{1}{h}.\]
\end{assumption}
\begin{assumption}
\label{ass:lip}
There exists non-negative finite constants $\beta>0$ and $0<L<+\infty$ such that:
\[\norm{\xx - \xx_t^*} \le \beta \Rightarrow\\ \norm{\nabla^2 f_t(\xx)-\nabla^2 f_t(\xx_t^*)}\le L\norm{\xx-\xx_t^*}.\] 
\end{assumption}
\begin{assumption}
\label{ass:lip2}
There exists $0<l<+\infty$ such that:
\[\norm{f_t(\xx)-f_t(\xx^*_t)}\le l\norm{\xx-\xx^*_t}.\]
\end{assumption}

Assumption \ref{ass:inv} imposes an upper bound on the norm of the inverse Hessian at the optimum. This implies that the Hessian's eigenvalues can be positive or negative but must be bounded away from zero. {For convex loss functions, this is equivalent to strong convexity which is a common assumption in OCO \cite{Shalev, OnlineNewton, OnlineEpsilonProof}. } Assumptions \ref{ass:lip} and \ref{ass:lip2} are local Lipschitz continuity conditions on the objective function and its Hessian around the optimum.

\subsection{Reduced function identities}

We now provide two lemmas which characterize the reduced function $\tilde{f}_t$.

\begin{lemma}\label{lemz}\cite[Section 10.2.3]{boyd}
Suppose
\begin{enumerate}
    \item $\exists \FF_t \text{ such that } \mathcal{N}(\AA_t) = \mathcal{R}(\FF_t);$
    \item $\exists \xh \text{ such that } \AA_t\xh=\bb_t;$
    \item $\AA_t\xx_t-\bb_t=\mathbf{0}.$
\end{enumerate}
Consider the Newton step applied to the reduced function $\tilde{f}_t(\zz_t)$:
$\Delta \zz_t = -\nabla^2\tilde{f_t}(\zz_t)^{-1}\nabla \tilde{f_t}(\zz_t)$.
Then the following identity holds :
\begin{align}
    \Delta \xx_t = \FF_t \Delta \zz_t.
\end{align}
\end{lemma}
Lemma \ref{lemz} implies that the Newton step applied to the constrained problem coincides with the Newton step applied to the reduced problem. By setting $\FF_t=\overline{\FF}_t$, we obtain $\norm{\Delta\xx_t} = \norm{\Delta\zz_t}$.

The second lemma characterizes the local strong convexity and Lipschitz continuity of the reduced function.
\begin{lemma}
\label{lem:const}
Suppose Assumptions \ref{ass:inv} and \ref{ass:lip} hold. Then we have:
\begin{align}
\label{constants1}
    &\norm{\nabla^2 \tilde{f}_t(\zz_t^*)^{-1}} \le \frac{1}{\sigma_{\min}(\FF_t)^2h}\\
    &\norm{\zz-\zz^*_t}\le \frac{\beta}{\norm{\FF_t}}\Rightarrow\nonumber\\
    \label{constants2}
    &\quad\norm{\nabla^2\tilde{f_t}(\zz)-\nabla^2\tilde{f_t}(\zz_t^*)}\le L\norm{\FF_t}^3\norm{\zz-\zz_t^*},
\end{align}
\end{lemma}
where $\sigma_{\min}(\FFt)$ is the minimum singular value of $\FFt$.

\noindent\begin{proof} 
Differentiating $\tilde{f}_t$ twice, we obtain, 
\[\nabla^2\tilde{f_t}(\zz) = \FF_t\nabla^2f_t(\xx_t^*)\FF_t^\top.\] 
The inverse Hessian of $\tilde{f}_t$ is thus upper bounded by:
\begin{align*}
    \norm{\nabla^2\tilde{f}_t(\zz^*)^{-1}} &= \norm{\FFt(\FFt^\top\FFt)^{-1}\nabla^2f_t(\xx_t^*)^{-1}(\FFt^\top\FFt)^{-1}\FFt^\top}\\
    &\le \norm{(\FFt^\top\FFt)^{-1}\FFt^\top}^2\norm{\nabla^2f_t(\xx_t^*)^{-1}}
    \\
    &\le \frac{1}{\sigma_{\min}(\FFt)^2h},
\end{align*}
which is \eqref{constants1}.\\
For \eqref{constants2}, we have,
\begin{align*}
    \norm{\nabla^2\tilde{f_t}(\zz)-\nabla^2\tilde{f_t}(\zz_t^*)} &= \big\|\FFt\big(\nabla^2f_t(\FFt\zz+\xh)-\\&\quad\quad\quad\nabla^2f_t(\FFt\zz^*_t+\xh)\big)\FFt^{\top}\big\|\\
    &\le \norm{\FFt}^2\big\|\nabla^2f_t(\FFt\zz+\xh)-\\&\quad\quad\quad\quad\nabla^2f_t(\FFt\zz^*_t+\xh)\big\|\\
    \quad\quad\quad\quad\quad&\le \norm{\FFt}^2 L\norm{\FFt(\zz-\zz^*_t)}\\
    &\le L\norm{\FFt}^3\norm{\zz-\zz^*_t},
\end{align*}
where the last inequalities follow from the Lipschitz continuity of $\nabla^2f_t$ and the definition of $\zz_t$, respectively.
\end{proof}

\subsection{Feasible Newton update}

Using Lemmas \ref{lemz} and \ref{lem:const}, we derive the following lemma for the original constrained problem:
\begin{lemma}[Equality-constrained Newton identities]
\label{lem:ECN}
Suppose Assumptions \ref{ass:inv}, \ref{ass:lip} hold and:
\begin{enumerate}
    \item $\norm{\xx_t-\xx_t^*}\le \min \left\{\beta, \frac{h}{2L}\right\};$
    \item $\AA_t\xx_t-\bb_t=\mathbf{0}.$
\end{enumerate}
    Then we have the following two identities for \texttt{OEN}:
    \begin{align}
    \label{id:1}
    \norm{\xx_{t+1}-\xx_t^*} &< \norm{\xx_t-\xx_t^*};\\
    \label{id:2}
    \norm{\xx_{t+1}-\xx_t^*} &\le \frac{2L}{h}\norm{\xx_t-\xx_t^*}.
\end{align}
\end{lemma}
{The first inequality guarantees the next iterate is strictly closer to the optimum compared to the current iterate. The second inequality provides an upper bound on this value.}

\begin{proof}
By the definition of the \texttt{OEN} update and Lemma~\ref{lemz} we have:
\[\xx_{t+1}-\xx_t^* = \xx_t-\xx_t^*-\FF_t\nabla^2\tilde{f}_t(\zz_t)^{-1}\nabla\tilde{f}_t(\zz_t).\]
Rearranging and letting $\FF_t=\FFo$, we have
\begin{align*}
    \xx_{t+1}-\xx_t^* &= \xx_t-\xx_t^*-\\
    &\quad\quad\quad\FFo\nabla^2\tilde{f}_t(\zz_t)^{-1}\big(\nabla\tilde{f}_t(\zz_t)-\nabla\tilde{f}_t(\zz_t^*)\big)\\
    &= \xx_t-\xx_t^*-\FFo\nabla^2\tilde{f}_t(\zz_t)^{-1}\cdot\\
    &\quad\quad\quad\int_0^1\nabla^2\tilde{f}_t\big(\zz_t+\tau(\zz_t^*-\zz_t)\big)(\zz_t^*-\zz_t)\text{d}\tau.
\end{align*}
From the symmetry of the Hessian and its inverse we have,
\begin{align*}
    \xx_{t+1}-\xx_t^* &= \xx_t-\xx_t^*-\FFo(\zz_t^*-\zz_t)\nabla^2\tilde{f}_t(\zz_t)^{-1}\cdot\\
    &\quad\quad\quad\quad\quad\int_0^1\nabla^2\tilde{f}_t\big(\zz_t+\tau(\zz_t^*-\zz_t)\big)\text{d}\tau\\
    &= (\xx_t^*-\xx_t)\nabla^2\tilde{f}_t(\zz)^{-1}\cdot\\
    &\quad\quad\int_0^1\Big[\nabla^2\tilde{f}_t\big(\zz+\tau(\zz^*-\zz)\big)-\nabla^2\tilde{f}_t(\zz)\Big]\text{d}\tau.
\end{align*}
Taking the norm on both sides and using Lemmas \ref{lemz} and \ref{lem:const} yields
\begin{align*}
    \norm{ \xx_{t+1}-\xx_t^*} &\le \norm{\xx_t^*-\xx_t}\norm{\nabla^2\tilde{f}_t(\zz)^{-1}}\mathbf{\cdot}\nonumber\\
    &\quad\int_0^1\norm{\nabla^2\tilde{f}_t\big(\zz+\tau(\zz^*-\zz)\big)+\nabla^2\tilde{f}_t(\zz)}\text{d}\tau\nonumber
    \end{align*}
    \begin{align}
    &\le \norm{\xx_t^*-\xx_t}\norm{\nabla^2\tilde{f}_t(\zz)^{-1}}\nonumber\\
    &\quad\quad\quad\int_0^1 2L\norm{\FFo}^2\tau\norm{\xx_t^*-\xx_t}\text{d}\tau\nonumber\\
    &\le \norm{\nabla^2\tilde{f}_t(\zz)^{-1}}L\norm{\xx_t^*-\xx_t}^2.
    \label{eq:1}
\end{align}
Using \cite[Lemma 2]{OnlineNewton}, we can bound the inverse Hessian as:
\begin{align}
    \norm{\nabla^2\tilde{f}(\zz)^{-1}} &\le \frac{1}{\sigma_{\min}(\FFo)^2h-L\norm{\FFo}^2\norm{\xx^*_t-\xx_t}}\nonumber\\
    &\le \frac{1}{h-L\norm{\xx^*_t-\xx_t}},
    \label{eq:2}
\end{align}
because $\sigma_{\min}(\FFo)=\norm{\FFo}=1$. Substituting (\ref{eq:1}) into (\ref{eq:2}) leads to:
\begin{align}
\label{eq:N}
    \norm{\xx_{t+1}-\xx_t^*} &\le \frac{L}{h-L\norm{\xx_t^*-\xx_t}}\norm{\xx_t^*-\xx_t}^2
\end{align}
Finally, (\ref{id:1}) and (\ref{id:2}) follows from \eqref{eq:N} and \cite[Lemma 2]{OnlineNewton}'s proof.
\end{proof}
\section{Online Equality-constrained Newton's Method}
We now present our online optimization methods for problems with time-independent and time-varying equality constraints. 
\subsection{Online Equality-constrained Newton's Method}
In this section, we propose the Online Equality-constrained Newton's Method (\texttt{OEN-M}) for online optimization subject to time-invariant linear equality constraints. This is the first online, second-order algorithm that admits constraints. \texttt{OEN-M} is presented in Algorithm~\ref{alg:OENM}.

\begin{algorithm}[h!]
\caption{Online Equality-constrained Newton's Method}\label{alg:OENM}
\begin{algorithmic}
\State\textbf{Parameters:} $\AA$, $\mathbf{b}$
\State\textbf{Initialization}: Receive $\xx_0 \in\mathbb{R}$ such that $\norm{\xx_0-\xx_0^*}\le\gamma=\min\left\{\beta,\frac{h}{2L}\right\}$ and $\AA\xx_0-\mathbf{b}=\mathbf{0}$
\FOR{$t=0,1,2...T$}
\State Play the decision $\xx_t$.
\State Observe the outcome $f_t(\xx_t)$.
\State Update decision:
\State   $\begin{bmatrix}\xx_{t+1}\\ \vv_t\end{bmatrix}$ $=\begin{bmatrix}{\xx}_{t}\\ \mathbf{0}\end{bmatrix}-\DD_t({\xx}_t)^{-1}\begin{bmatrix}\nabla f_t(\xx_t)\\\mathbf{0}\end{bmatrix}$.
\ENDFOR
\end{algorithmic}
\end{algorithm}

We now show that \texttt{OEN-M} has a dynamic regret bounded above by $O(V_T+1)$ and a null constraint violation.

\begin{theorem}
\label{regretOENM}
If Assumptions \ref{ass:inv}$-$\ref{ass:lip2} hold and the following conditions are respected:
\begin{enumerate}
    \item $\exists \xx_0 \text{ such that } \norm{\xx_0-\xx_0^*} \le \gamma = \min\{\beta, \frac{h}{2L}\}$;
    \item $\overline{v} \le \gamma - \frac{2L}{h}\gamma^2$.
\end{enumerate}
Then the regret $R_{\text{d}}(T)$ and the constraint violation $\text{Vio}(T)$ are bounded above by :
\begin{align}
\label{eq:regret}
    R_{\text{d}} (T) &\le \frac{lh}{h-2L\gamma}(V_T+\delta);\\
    \text{Vio}(T)&=0.
\end{align}
\end{theorem}

\begin{remark}
{The assumption that the decision-maker has access to $\xx_0$ such that $\norm{\xx_0-\xx_0^*}\le\gamma$ is standard in OCO \cite{OnlineConvexBandit, OnlineNewton, DynamicPolyak}. Essentially, this means obtaining a good starting estimate of the initial optimal solution is required. It is assumed that a good estimate can be obtained before the start of the online process from, for example, offline calculations or a previously implemented decision.}
\end{remark}

\begin{proof}
    Using Assumption \ref{ass:lip2}, the regret is bounded by:
\begin{equation}
\label{eq:reg1} R_{\text{d}}(T) =\sum\limits_{t=1}^T\big|f_t(\xx_t)-f_t(\xx_t^*)\big| \le l\sum\limits_{t=1}^T \norm {\xx_t-\xx_t^*}.\end{equation}

Rearranging (\ref{eq:reg1})'s sum we obtain:
\begin{align*}
    \sum\limits_{t=1}^T \norm {\xx_t-\xx_t^*} &=  \sum\limits_{t=1}^T \norm {\xx_t-\xx_{t-1}^*+\xx^*_{t-1}-\xx_t^*}\\
    &\le \sum\limits_{t=1}^T \norm {\xx_t-\xx_{t-1}^*} + \sum\limits_{t=1}^T \norm {\xx_t^*-\xx_{t-1}^*}\\
    &\le \sum\limits_{t=0}^{T-1} \frac{2L}{h}\norm{\xx_t-\xx_t^*}^2 + V_T\\
    &\le \sum\limits_{t=1}^{T} \frac{2L}{h}\gamma\norm{\xx_t-\xx_t^*} + V_T + \delta,
\end{align*}
where $\delta = \frac{2L}{h}\gamma(\norm{\xx_0-\xx_0^*}-\norm{\xx_T-\xx_T^*})$.  Solving for $\norm{\xx_t-\xx_t^*}$, we have:
\begin{align}
\label{eq:reg2}
    \sum\limits_{t=1}^T \norm {\xx_t-\xx_t^*} &\le \left(1-\frac{2L}{h}\gamma\right)^{-1}(V_T+\delta).
\end{align}
This implies that the dynamic regret is bounded above by:
\begin{align*}
    R_{\text{d}}(T) \le \frac{lh}{h-2L\gamma}(V_T+\delta),
\end{align*}
and hence, $R_{\text{d}}(T) \le O(V_T+1)$.

As for the constraint violation, we have that $\xx_0$ is feasible by assumption. Because every \texttt{OEN} update is such that $\AA\Delta\xx=0$, every subsequent decision will also be feasible. We thus have:
\begin{align*}
    \text{Vio}(T) &= \sum\limits_{t=1}^{T}\norm{\AA\xx_t-\bb}=0,
\end{align*}
which completes the proof.
\end{proof}

\subsection{Online Projected Newton's Method}
We now consider online optimization problems with time-varying equality constraints. \texttt{OEN} does not apply to this class of problems because the previously played decision might not be feasible under the new constraints. We propose the Online Projected Equality-constrained Newton's Method (\texttt{OPEN-M}) to address this limitation. \texttt{OPEN-M} consists of a projection of the previous decision onto the new feasible set followed by an \texttt{OEN} step from this point. \texttt{OPEN-M} is detailed in Algorithm \ref{alg:ProjNewton}.
\begin{algorithm}
\caption{Online Projected Eq.-const. Newton Method}\label{alg:ProjNewton}
\begin{algorithmic}
\State\textbf{Initialization}: Receive $\xx_0 \in\mathbb{R}$ such that $\norm{\xx_0-\xx_0^*}\le\gamma=\min\left\{\beta,\frac{h}{2L}\right\}$
\FOR{$t=0,1,2...T$}
\State Play the decision $\xx_t.$
\State Observe the outcome $f_t(\xx_t)$ and constraints $\AA_t, \bb_t$.
\State Project $\xx_t$ onto the feasible set: 
\State $\tilde{\xx}_t = \xx_t+\AA_t^\top(\AA_t\AA_t^\top)^{-1}(\bb_t-\AA_t\xx_t).$
\State Update decision:
\State   $\begin{bmatrix}\xx_{t+1}\\ \vv_t\end{bmatrix}$ $=\begin{bmatrix}\tilde{\xx}_{t}\\ \mathbf{0}\end{bmatrix}-\DD_t(\tilde{\xx}_t)^{-1}\begin{bmatrix}\nabla f_t(\tilde{\xx}_t)\\\mathbf{0}\end{bmatrix}.$
\ENDFOR
\end{algorithmic}
\end{algorithm}
We now analyse the performance of \texttt{OPEN-M}.
\begin{theorem} Suppose Assumptions \ref{ass:inv}$-$\ref{ass:lip2} hold and:
\begin{enumerate}
    \item $\exists \xx_0 \text{ such that } \norm{\xx_0-\xx_0^*} \le \gamma = \min\{\beta, \frac{h}{2L}\}$;
    \item $\overline{v} \le \gamma - \frac{2L}{h}\gamma^2$;
    \item $\exists a>0\text { such that }\norm{\AA_t}\le a \quad\forall t$.
\end{enumerate}
Then the regret $R_{\text{d}}(T)$ and the constraint violation $\text{Vio}(T)$ of \texttt{OPEN-M} is bounded above by :
\begin{align}
\label{eq:OPENreg}
    R_{\text{d}} (T) &\le \frac{lh}{h-2L\gamma}(V_T+\delta)\\
    \label{eq:OPENcons}
    \text{Vio}(T)&\le \frac{ah}{h-2L\gamma}(V_T+\delta).\\\nonumber
\end{align}
\end{theorem}
\begin{proof}
    We first show that the following inequality holds:
\begin{equation}\label{eq:proj}\norm{\tilde{\xx}_t-\xx_t^*}\le \norm{\xx_t-\xx_t^*}.\end{equation}
Since $\tilde{\xx}_t$ is the projection of $\xx_t$ onto the feasible set at time $t$, $\xx_t-\xx_t^*$ and $\tilde{\xx}_t-\xx_t^*$ are orthogonal. It follows that:
$\norm{\tilde{\xx}_t-\xx_t^*}^2 = \norm{\xx_t-\xx_t^*}^2-\norm{\tilde{\xx}_t-\xx_t}^2.$
Because $\norm{\tilde{\xx}_t-\xx_t}\ge 0$, we have $\norm{\tilde{\xx}_t-\xx_t^*} \le \norm{\xx_t-\xx_t^*}$, which is~\eqref{eq:proj}.

This implies that $\norm{\tilde{\xx}_t-\xx_t^*} \le \gamma$ which means the projected decision $\tilde{\xx}_t$ satisfies all the requirements for \texttt{OEN-M}. The same analysis as for \texttt{OEN-M} therefore holds for \texttt{OPEN-M} and the same regret bound is obtained, thus leading to \eqref{eq:OPENreg}.

As for the constraint violation, we recall:
\[\text{Vio}(T) = \sum\limits_{t=1}^{T}\norm{\AA_t\xx_t-\bb_t}.\]
Using the fact that $\xx_t^*$ is feasible for every timestep,
\begin{align*}
    \text{Vio}(T) &= \sum\limits_{t=1}^{T}\norm{\AA_t\xx_t-\mathbf{b}_t-(\AA_t\xx_t^*-\bb_t)}\\
     &= \sum\limits_{t=1}^{T}\norm{\AA_t(\xx_t-\xx_t^*)}.
\end{align*}
By the Cauchy-Swartz inequality and the bound on $\AA_t$:
\begin{align*}
    \text{Vio}(T) &\le \sum\limits_{t=1}^{T}\norm{\mathbf{A}_t}\norm{\xx_t-\xx_t^*}\\
    &\le a\sum\limits_{t=1}^{T}\norm{\xx_t-\xx_t^*}\\
    &\le \frac{ah}{h-2L\gamma}(V_T+\delta),
\end{align*}
where the last inequality follows from (\ref{eq:reg2}). This yields (\ref{eq:OPENcons}) and completes the proof.
\end{proof}

\begin{remark}Because a closed-form projection step is possible for \texttt{OPEN-M}, the algorithmic time-complexity is the same as for \texttt{OEN-M}. This is because the time-complexity of \texttt{OEN-M} and \texttt{OPEN-M} are dominated by the matrix inversion step which is $O(n^5\log(n))$ in the general case. Thus, there is no additional burden to \texttt{OPEN-M}'s projection step which is $O(n^3)$ in the worst case.
\end{remark}

{\begin{remark}\label{remark2}
\texttt{OPEN-M} possesses the tightest dynamic regret bounds of any previously proposed online equality-constrained algorithm in the literature \cite{MOSP, QueueConstrained,LagrangianTime-Varying, DynamicPolyak}. Under the standard assumption in OCO that the variation of optima $V_T$ is sublinear \cite{LagrangianTime-Varying}, constraint violation and regret will be sublinear. The method is also parameter-free which eliminates the need for time-dependent step sizes and hyperparameter tuning, e.g., the step size in gradient-based methods. The time horizon during which the algorithm is used is also arbitrary and does not need to be defined before execution. These advantages provide ample justification for the additional complexity of the inversion step. \end{remark} }

\section{Numerical Experiment}
We now illustrate the performance of \texttt{OPEN-M} and use it in an optimal network flow problem. This type of problem can model electric distribution grids when line losses are considered as negligible \cite{NetworkFlow, NetFlow2}. In this context, a convex, quadratic cost is most commonly used. Note that, \texttt{OPEN-M} is also applicable to non-convex cost functions.

Consider the network flow problem over a directed graph $\mathcal{G}=(\mathcal{M}, \mathcal{L})$ with nodes $\mathcal{M}$ and directed edges $\mathcal{L}$. At every timestep $t$, load (sink) nodes, $i\in \mathcal{D}$, require a power supply $b_t^{i}$ and generator (source) nodes, $j\in \mathcal{P}$, can produce a positive quantity of power. The power is distributed through the edges of the graph. The decision variable is $\xx_t$ and models the power flowing through each edge. Assuming no active power losses, the power balance at each node leads to the constraint:
$\AA\xx_t=\bb_t$ where $b_t^{(i)}$ is the power demand at each load node and:
\begin{align*}\AA_{(l, i)} = \begin{cases}1, &\text{ if edge } l \text{ enters node } i\\-1, &\text{ if edge } l \text{ leaves node } i\\0, &\text{ else}.\end{cases}\end{align*}

A numerical example is provided next using a fixed, radial network composed of 15 nodes connected via 30 arcs. A single power source is located at the root of the network. Every node's load is chosen independently as: $b_t^i=\frac{\zeta}{\sqrt{t}}+10$ where $t$ denotes the round and $\zeta$ is uniformly sampled in $[0,5]$. The cost function $f_i$ for each arc $i$ is convex and adopts the following form: $f_i(x)=\alpha_i \mathrm{e}^{\beta_i |x|}$. The parameters $\alpha_i$ and $\beta_i$ are chosen following: $\alpha_i = \frac{\eta}{\sqrt{t}}+1$ and $\beta_i=\frac{\gamma}{\sqrt{t}}+2$ where $t$ denotes the round and $\eta$ and $\gamma$ are uniformly sampled in $[0,10]$ at every round. This loss function is chosen because it is harder to solve than a quadratic function and yet approximately models electric grid costs. The temporal dependence of the parameters ensures that the total variation of optima ($V_T$) is bounded and sublinear. The time horizon is set to $T=2500$. The fixed nature of the network and the diagonal Hessian matrix means that the inversion step only has to be done once. Note that \texttt{OPEN-M} also admits time-varying network topologies, i.e., using $\AA_t$ instead of $\AA$.

\begin{figure}
    \centering
        \includegraphics[width=0.45\textwidth]{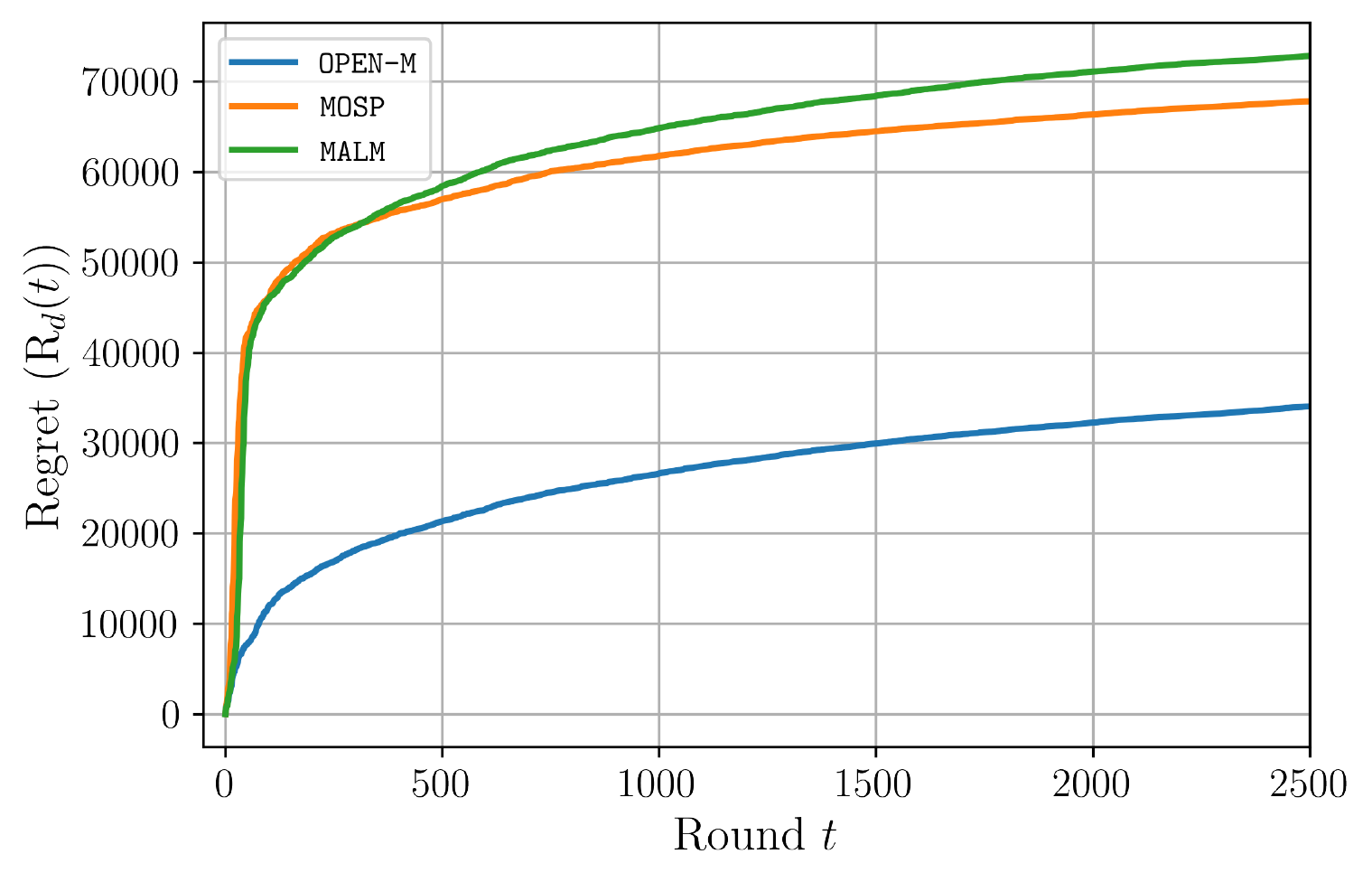}
        \vspace{-0.5cm}
        \caption{Dynamic regret comparison of \texttt{OPEN-M}, \texttt{MOSP} and \texttt{MALM}}
        \label{fig:regret}
\end{figure}

\begin{figure}
    \centering
    \includegraphics[width=0.45\textwidth]{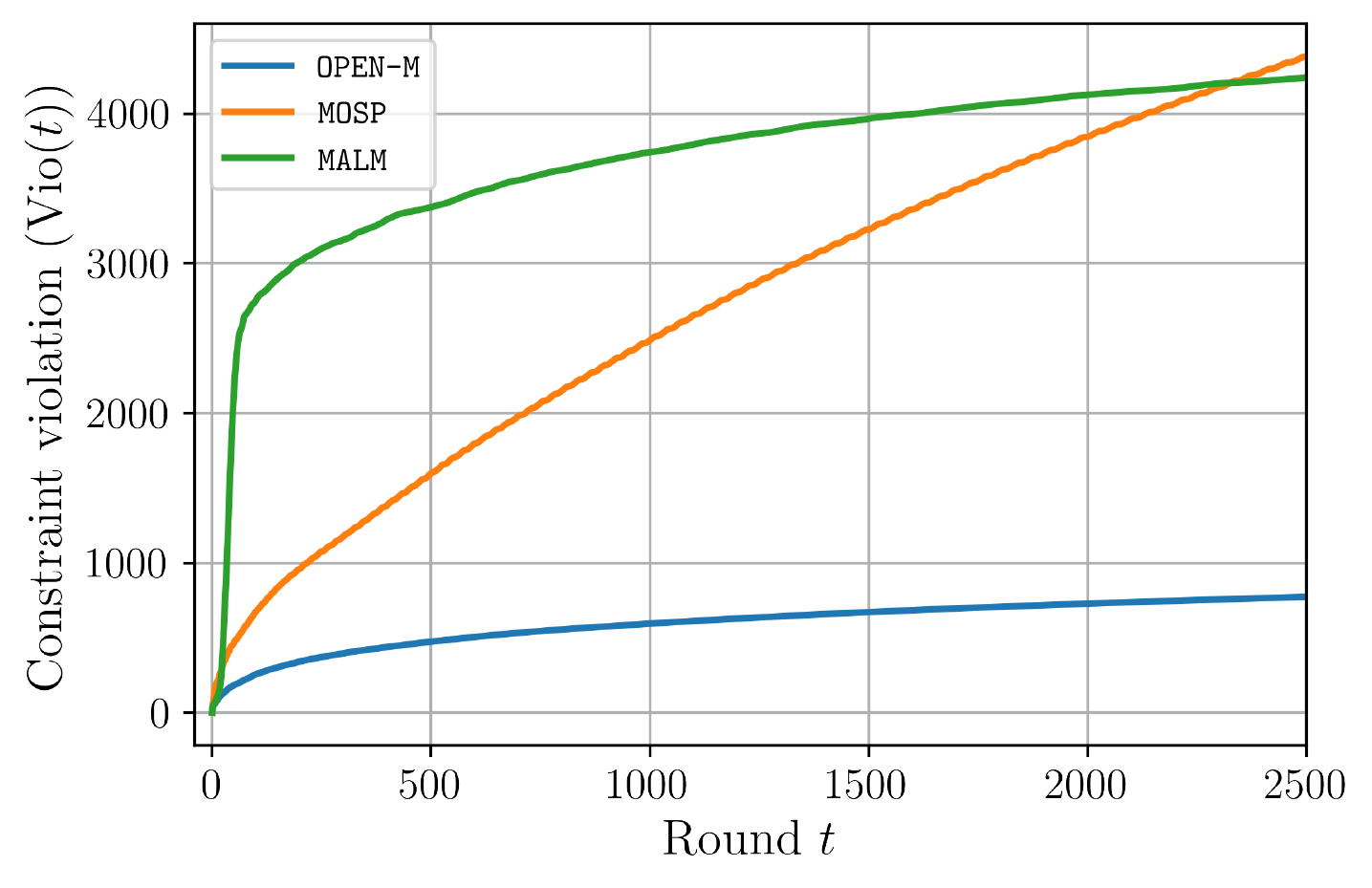}
    \vspace{-0.5cm}
    \caption{Constraint violation comparison of \texttt{OPEN-M}, \texttt{MOSP} and \texttt{MALM}}
    \label{fig:const}
    \vspace{-0.5cm}
\end{figure}

We use \texttt{MOSP} from \cite{MOSP} and the model-based augmented Lagrangian method (\texttt{MALM}) from \cite{LagrangianTime-Varying} as benchmarks to establish \texttt{OPEN-M}'s performance. Because these algorithms can only admit inequality constraints, the equality constraint $\AA\xx_t=\bb_t$ is relaxed to $\AA\xx_t-\bb_t\le0$. This ensures that there is enough power at every node but lets the operator over-serve loads. This relaxation is mild because the constraint should be active at the optimum given that costs are minimized.
{Dynamic regret, defined in \eqref{RegretDef}, and constraint violation, defined in \eqref{VioDef}, for this problem are presented in Figures \ref{fig:regret} and \ref{fig:const}, respectively.}

We observe sublinear dynamic regret and constraint violation from all three algorithms illustrating they are well-adapted to this problem. We remark that \texttt{OPEN-M} exhibits a lower regret than both the \texttt{MOSP} and \texttt{MALM} algorithms. Indeed, the dynamic regret of \texttt{OPEN-M} is an order of magnitude smaller. \texttt{OPEN-M} also has significantly better constraint violation compared to the other two algorithms.

\section{Conclusion}

In this paper, a second-order approach for online constrained optimization is developed. Under linear time-varying equality constraints, the resulting algorithm, \texttt{OPEN-M}, achieves simultaneous $O(V_T+1)$ dynamic regret and constraint violation bounds. These bounds are the tightest yet presented in the literature. A numerical network flow example is presented to showcase the performance of \texttt{OPEN-M} compared to other methods from the literature.

Considering the prevalence of interior-point methods in the offline optimization literature, an extension of the equality-constrained Newton's method which admits inequality constraints \cite{renegar}, a similar extension can be envisioned for \texttt{OPEN-M}. A second-order approach to an online optimization problem with time-varying inequalities has the potential to improve current dynamic regret and constraint violation bounds.


\bibliographystyle{IEEEtran}
\bibliography{Ref}

\end{document}